\documentclass[english]{IWMS}

\newcommand{\ba}{\begin{array}}
\newcommand{\ea}{\end{array}}
\newcommand{\bt}{\begin{tabular}}
\newcommand{\et}{\end{tabular}}
\newcommand{\btb}{\begin{table}}
\newcommand{\etb}{\end{table}}
\newcommand{\bc}{\begin{center}}

\newcommand{\bea}{\begin{eqnarray}}
\newcommand{\eea}{\end{eqnarray}}
\newcommand{\Bea}{\begin{eqnarray*}}
\newcommand{\Eea}{\end{eqnarray*}}
\newcommand{\beq}{\begin{equation}}
\newcommand{\eeq}{\end{equation}}

\begin{document}


\begin{englishtitle}

\title{Antieigenvalue Analysis, New Applications: Continuum Mechanics, Economics, Number Theory}{} 

\author[1]{Gustafson, Karl}{Copyright \copyright 2015 by Karl Gustafson
All rights Reserved}



\address[1]{University of Colorado at Boulder, USA}

\maketitle

\begin{abstract}My recent book Antieigenvalue Analysis, World-Scientific, 2012, presented the theory of antieigenvalues from its inception in 1966 up to 2010, and its applications within those forty-five years to Numerical Analysis, Wavelets, Statistics, Quantum Mechanics, Finance, and Optimization. Here I am able to offer three further areas of application: Continuum Mechanics, Economics, and Number Theory.
\end{abstract}

\end{englishtitle}

\makecollection


\section{Introduction}

Antieigenvalue analysis \cite{1} is an operator trigonometry concerned with those vectors, called antieigenvectors, which are most-turned by a matrix or a linear operator A.  This is in contrast to the conventional eigenvalue analysis, which is concerned with those vectors, called eigenvectors, which are not turned at all by A.  Antieigenvalue theory may be usefully thought of as a variational theory, extending the variational Rayleigh-Ritz theory which characterizes eigenvectors, to an enlarged theory also characterizing antieigenvectors. 

Two key entities in the antieigenvalue theory are the first
antieigenvalue \beq \mu_{1} =  \cos \phi (A) =\min_{x \neq 0}\frac{
<Ax,x>}{ \parallel Ax \parallel \parallel x \parallel }\label{1.1}
\eeq 
and the related convex minimum
\beq
\nu_{1} = \sin \phi(A) =\min_{ \epsilon > 0} \parallel \epsilon A-I \parallel \mathrm{.}\label{1.2} 
\eeq
Here I will specialize the antieigenvalue theory to A, an $n \times n$ symmetric positive definite matrix. One generally has the fundamental relation 
\beq
\cos^{2} \phi (A) + \sin^{2}  \phi (A)=1 \mathrm{.}\label{1.3} 
\eeq
There are two maximally turned first antieigenvectors 
\beq
x_{\pm} = \left( \frac{\lambda _{n}}{\lambda_{1} + \lambda_{n}}\right) ^{\frac{1}{2}} x_{1} \pm \left( \frac{\lambda _{1}}{\lambda_{1} + \lambda_{n}}\right) ^{\frac{1}{2}} x_{n} \label{1.4} 
\eeq
where $0 < \lambda_{1} \leqq \lambda_{2} \leqq \ldots \leqq \lambda_{n}$ are the eigenvalues of A, and where $x_{1}$ is any norm-one eigenvector from the $\lambda_{1}$-eigenspace and $x_{n}$ is any norm-one eigenvector from the $\lambda_{n}$-eigenspace.  The antieigenvectors $x_{\pm}$ in (\ref{1.4}) have also been normalized to be of norm-one. For $n \times n$ symmetric positive definite A the expressions in (\ref{1.1}) and (\ref{1.2}) have useful explicit valuations as
\beq
\mu_{1} =\frac{2\sqrt{\lambda_{1}\lambda_{n}}}{\lambda_{1}+\lambda_{n}} \ ,\   \nu_{1}=\frac{\lambda_{n}-\lambda_{1}}{\lambda_{n}-\lambda_{1}}\mathrm{.}\label{1.5}
\eeq 



For further elaboration of the general antieigenvalue theory I refer to \cite{1}.  In particular, just as one may move up the successive eigenvalue ladder via the Rayleigh-Ritz variational quotient minimizations, one can analogously move up an antieigenvalue ladder $\mu_{2},\mu_{3}, ...$ via the variational quotient minimizations of (\ref{1.1}), thereby arriving at a decreasing sequence of critical turning angles $\phi_{k}(A)$.  But we won't need those higher antieigenvalues in the discussions of this paper. 

In Sections 2, 3 and 4, respectively, I will briefly summarize three new applications of the antieigenvalue analysis. A full development of each will appear in the three forthcoming papers (\cite{2},\cite{3},\cite{4}), respectively.

I would like to express my thanks to Jeff Hunter and Simo Puntanen as chairs of the IWMS-2015 for inviting me and to ILAS for designating me as their lecturer for the conference.
\\\\
\section{Continuum Mechanics(Granular Materials)}
The following extracts a main result from the forthcoming paper \cite{2}. My investigations there lean heavily on the recent paper \cite{5} which I would urge the reader to consult for further background on the modeling of granular materials. In \cite{2} I also go further to make connections to sandpile theories which have been useful to model self-organizing criticality in Statistical Mechanics \cite{6}. 

The paper \cite{5} explores the notion of (maximum) angle of repose for granular materials. On the other hand, my theory of antieigenvalues \cite{1} has as one of its essential ingredients the notion of (maximum) operator turning angle. Here is how to connect the two theories.   

Following \cite{5}, the equilibrium equations for a granular pile of local slope $\theta$ are
\beq \begin{array}{ccccc}
\partial_{x} \sigma_{xx} & + & \partial_{z} \sigma_{xz} & = & \rho g \sin\theta \\
\partial_{x} \sigma_{xx} & + & \partial_{z} \sigma_{xz} & = & \rho g \cos\theta
\end{array}
\mathrm{.}\label{2.1} 
\eeq
The stress tensor $\Sigma$ can be written in singular value decomposition\\
\beq
\Sigma = \left[ \begin{array}{cc}
\sigma_{xx} & \sigma_{xz}\\
\sigma_{xz} & \sigma_{zz}
\end{array} \right] =
\left[ \begin{array}{cc}
\cos\psi & -\sin\psi\\
\sin\psi & \cos\psi
\end{array} \right]
\left[ \begin{array}{cc}
\sigma_{1} & 0\\
0 & \sigma_{2}
\end{array} \right]
\left[ \begin{array}{cc}
\cos\psi & \sin\psi\\
-\sin\psi & \cos\psi
\end{array} \right]
\eeq \label{2.2}
where $\sigma_{1} \geqq \sigma_{2} > 0$ are the principal stresses and where $\psi$ gives the principal directions.\\
By considering a plane within the material and the normal and tangential stresses upon it in terms of the coefficient of friction of the material and a corresponding angle $\delta$ of internal friction, it is deduced in \cite{5} that the largest sustainable angle of repose $\theta$ is given by 
\beq
\sin\theta = \frac{\tau}{\sigma}, where \  \sigma=\dfrac{\sigma_{1}+\sigma_{2}}{2}, \tau =\dfrac{\sigma_{1}-\sigma_{2}}{2}\mathrm{.}\label{2.3}\eeq 
Some assumptions in this modeling of the discrete by the continuous have of course been made. Among those are a linear dependence on the vertical z direction and a stress-free condition at the pile's surface at $z=0$. 

In \cite{2} I am able to see and recast this continuum mechanical model for stable granular material piles into my antieigenvalue theory.  The key is to remember that the unique $\epsilon$ for which the minimum in (\ref{1.2}) is attained in known \cite{1} to be $\epsilon_{m} = \frac{2}{\sigma_{1}+\sigma_{2}}$ for a $2 \times 2$ matrix with singular values $\sigma_{1}$ and $\sigma_{2}$.  Then straightforward calculations confirm that 
\beq
\sin\phi(\Sigma )= \Vert \epsilon_{m}\Sigma - I \Vert = \frac{\tau}{\sigma} \mathrm{.}\label{2.4} 
\eeq
Theorem (\cite{2}, Thm 3.1). The (maximum) angle of repose $\theta$, the material angle of friction $\delta$, and the (maximum) turning angle $\phi(\Sigma)$ of the stress tensor, are the same.\\\\ 
\section{Economics (Capital Asset Pricing Model)}
The following extracts from the forthcoming paper \cite{3}.  I presented a preliminary version at the IWMS2013 conference in Toronto. The paper \cite{7} ensuing from that conference was accepted and may be accessed online although it is not in print yet at this writing.  It contains a rather complete articulation, in order, of randomness, risk, and reward in financial markets, along with a number of important bibliographical references.  I refer you here to the excellent book \cite{8} for aspects of the CAPM as it is used in portfolio design theory and in high frequency trading, topics I have followed now for a considerable number of years.

The main finding which I exposed at the IWMS2013 and in the paper \cite{7} and a bit earlier in the book [\cite{1}, p.182] is that the Sharpe Ratio, which is a key tool used in the CAPM theory and more importantly in practice, can be related to my first antieigenvalue $\mu_{1}$. I first observed this in 1994 when I entered the financial engineering mathematical world during a very small consulting task, but did not work out the details until \cite{7}.  Here, in brief, is how that relationship may be seen.

The Capital Asset Pricing model assumes the Efficient Market
Hypothesis and then tells you to measure the return-to-risk of your portfolio against the market.  From the assumption that the full market has optimized the return-to-risk, your Sharpe ratio 
\beq
S=\frac{E[r]}{\sigma[r]},\label{3.1} 
\eeq where $E[r]$ is the average return over a number of periods and $\sigma[r]$ is the corresponding standard deviation, will not be greater than that of the whole (e.g, think indexing) market's Sharpe ratio.  See especially [\cite{8}, fig 5.2, p 55], to picture Sharpe ratios as mean-variance slopes.  Here I am just dropping the risk-free return rate $R_{f}$ from numerators, as it is effectively zero these days anyway.

Suppose now we look at the last two years of annualized returns $r_{1}$ and $r_{2}$.  We may form the usual (arithmetic) Sharpe ratio $S_{AM}= \frac{r_{1}+r_{2}}{2\sigma}$ and also a (geometric) Sharpe ratio $S_{GM}=\frac{\sqrt{r_{1}r_{2}}}{\sigma}$ and upon dividing the latter by the former we arrive at 
\beq
G=\frac{S_{GM}}{S_{AM}}= \frac{2\sqrt{r_{1}r_{2}}}{r_{1}+r_{2}} \label{3.2}
\eeq 
which is my first antieigenvalue $\mu_{1}$ as seen from (\ref{1.5}).  Further details and refinements may be found in \cite{3,7}. There I also begin an accompanying treatment of geometric versus arithmetic portfolio design.  A referee of \cite{7} also suggested possible connections to the currently important financial economic issues concerning realized volatilities, and I am in the process of such investigations in the paper under preparation \cite{3}.

\section{Number Theory (Pythagorean Triples)}

Here is a very interesting and new explicit connection of my antieigenvalue Theory and its operator trigonometry, to number theory.  More details will be given in a paper under preparation \cite{4}, where a number of important further ramifications will also be developed.


Let me first point out and emphasize that when I
first originated the antieigenvalue theory almost fifty years ago, I was coming from semi-group perturbation theory which had led me to a question of when an operator product BA would remain (real) positive, given positive A under multiplicative perturbation by positive B.  For rather general semi-group generators A, and bounded B, I found the operator theoretic sufficient condition
\beq
\sin\phi(B) \leqq \cos\phi(A)\mathrm{.}\label{4.1} 
\eeq
Then by using variational techniques on the expression (\ref{1.1}) in conjunction with convexity techniques on the expression (\ref{1.2}) I found the explicit valuations (\ref{1.5}) for $n \times n$ symmetric positive definite matrices.  See \cite{1} for more details and history. 

Now the new connection to Number Theory, which I only recently discovered.  Given two arbitrary relatively prime positive integers $m$ and $n$, with $m > n$, one of them being even, the other odd, then the numbers   
\beq
a=2mn,\ b=m^{2}-n^{2},\ c=m^{2}+n^{2} \label{4.2}
\eeq
form a primitive Pythagorean triple:
\beq
a^{2}+b^{2}=c^{2}\mathrm{.}\label{4.3} 
\eeq
This sufficient condition is also necessary.  For more details see \cite{9}.  This construction and characterization of Pythagorean triples is often called Euclid's Formula. 

I may now form the matrix (and its similarity class of matrices with the same eigenvalues $m^{2}$ and $n^{2})$
\begin{equation}
A= \begin{bmatrix}
n^{2} & 0\\
0 & m^{2}
\end{bmatrix}\mathrm{.}\label{4.4}
\end{equation} 
Immediately from my matrix operator trigonometry \cite{1} and the expressions (\ref{1.5}) we have
\beq
\cos\phi(A) =\frac{2mn}{m^{2}+n^{2}},\ \sin\phi(A)= \frac{m^{2}-n^{2}}{m^{2}+n^{2}}\mathrm{.} \label{4.5} 
\eeq
Proposition (\cite{4}).  Euclid's Formula for Pythagorean Triples is a special case of my operator trigonometry.  

We may propose to call these matrices $ A_{m,n} $, Pythagorean Triple Matrices.  Their maximum turning angles may be called special Pythagorean turning angles $ \phi_{m,n}(A) $.  Their corresponding normalized Pythagorean antieigenvectors are 
\beq
x_{\pm} =\left( \frac{m^{2}}{m^{2}+n^{2}}\right)^{\frac{1}{2}}\left[ \begin{array}{c}
1\\
0
\end{array}\right] \pm \left( \frac{n^{2}}{m^{2}+n^{2}}\right)^{\frac{1}{2}}\left[ \begin{array}{c} 
0\\
1
\end{array}\right]= \frac{1}{\sqrt{m^{2}+n^{2}}} \left[ \begin{array}{c}
m\\
\pm n
\end{array}\right]\mathrm{.}\label{4.6} 
\eeq
We know of course there are an infinite number of these Pythagorean angles, which are now embedded within my antieigenvalue operator trigonometry. 

This new connection of my antieigenvalue analysis to the Pythagorean triple number theory may be seen to have other interesting manifestations.  Here is another one, couched in the terminology of algebraic geometry.  Let $ x=\left( \frac{m}{n}, 0 \right)  $ be a point on the x-axis.  It's stereographic projection onto the unit circle becomes, now seen operator-theoretically, the point
\beq P= \left( \frac{2mn}{m^{2}+n^{2}},\
\frac{m^{2}-n^{2}}{m^{2}+n^{2}} \right) = \left(  \cos\phi
(A_{m,n}),\ \sin\phi (A_{m,n}) \right)\mathrm{.}\label{4.7} \eeq
The stereographic point of view comes from a treatment of spinors and twistors \cite{10}.
\\

\end{document}